\documentclass[14pt]{extarticle}
\usepackage{amsfonts,amssymb,eucal}
\usepackage{mathenv}
\usepackage{amsmath}
\usepackage{mathtext}
\usepackage[T2A]{fontenc}
\usepackage[cp1251]{inputenc}
\usepackage[russian]{babel}
\usepackage[dvips]{graphicx,rotating}
\usepackage{amssymb}
\usepackage{amsthm}
\usepackage{bm}
\usepackage{latexsym}
\usepackage{titlesec}
\usepackage{amsmath}
\usepackage{amsfonts}
\usepackage{cite}
\usepackage{eufrak}
\usepackage{eucal}
\usepackage{color}
\usepackage[a4paper, left=30mm, right=10mm, bottom=20mm, head=0mm]{geometry}
\begin{document}
 \thispagestyle{empty}

\begin{flushleft}
УДК 517.986.62
\end{flushleft}

\begin{center}
 ФУНКЦИИ ОГРАНИЧЕННОЙ СРЕДНЕЙ ОСЦИЛЛЯЦИИ И ГАНКЕЛЕВЫ
ОПЕРАТОРЫ НА КОМПАКТНЫХ АБЕЛЕВЫХ ГРУППАХ\footnote{Работа выполнена при поддержке второго из авторов Государственной программой научных исследований Республики Беларусь, проект № 20111164.}
\end{center}

\begin{center}
 {\bf  Р. В. Дыба, А. Р. Миротин}//
 amirotin@yandex.ru
 \end{center}

{\footnotesize Рассматривается  обобщение понятия функции ограниченной средней осцилляции и ганкелева оператора на случай компактных абелевых групп с линейно упорядоченной  группой характеров.  Дается описание пространств функций ограниченной средней осцилляции и  функций ограниченной средней осцилляции аналитического типа на таких группах в терминах ограниченности соответствующих операторов Ганкеля в предположении, что группа характеров содержит наименьший положительный элемент.}

{\footnotesize Ключевые слова:  оператор Ганкеля, ограниченный оператор,  ограниченная средняя осцилляция, линейно упорядоченная абелева группа, компактная абелева группа}.

 {\footnotesize R. V. Dyba, A. R. Mirotin. Functions of  bounded mean oscillation and Hankel operators  on compact abelian groups.

 Generalization of  functions of bounded mean oscillation and Hankel operators to the case of compact abelian groups with linearly ordered dual is considered. Spaces of functions of bounded mean oscillation and of  bounded mean oscillation of analytic type on such groups are described in terms  of boundedness of corresponding Hankel operators under the assumption that the dual group contains a minimal positive element.}


\begin{center}
{\bf 1. Введение}
\end{center}

Пространство $BMO(G)$ функций ограниченной средней осцилляции на
компактной абелевой группе $G$  были определены в работе \cite{DM} (относительно классического
случая группы вращений окружности см., например, \cite{Gar}), а ганкелевы операторы на таких группах,
рассматриваемые ниже, --- в \cite{Dyba}; другая версия рассматривалась в
\cite{YCG}  (по поводу классической теории операторов Ганкеля см. \cite{Nik1}, \cite{Pel}). Теория этих операторов тесно связана с обобщениями
на группы операторов Тёплица (см., например,\cite{SbMath} и литературу там), а также операторов Винера-Хопфа
\cite{Adukov}. Данная работа посвящена изучению  связи теории пространств функций ограниченной средней осцилляции и   функций ограниченной средней осцилляции аналитического типа на  группе $G$ с теорией  операторов Ганкеля в пространствах Харди
 $H^2(G)$. Основной результат (теорема 1)  дает  описание пространств $BMO(G)$  и $BMOA(G)$  в терминах ограниченности соответствующих операторов Ганкеля при условии, что группа характеров группы $G$ содержит наименьший положительный элемент (класс
 таких групп весьма широк и описан в  \cite{Adukov}, лемма 3.2). Остальные утверждения работы носят вспомогательный характер. Обобщение на рассматриваемый случай теоремы двойственности Феффермана
 будет дано в другой работе авторов. При этом будет выполнена программа,
 для классического случая намеченная в \cite{Nik1}, с. 189 -- 190.

\begin{center}
{\bf 2. Обозначения и вспомогательные сведения}
\end{center}

Всюду ниже   $G$ есть нетривиальная связная компактная абелева
группа с  нормированной мерой Хаара $dx$ и линейно упорядоченной
группой характеров $X$ \cite{Pont}, $X_+$ --- положительный конус в  $X$.
Другими словами,  $X$ есть дискретная абелева группа, в которой выделена подполугруппа $X_+$,
содержащая единичный характер $\textbf{1}$ и такая, что $X_+\cap
X_+^{-1}=\{\textbf{1}\}$ и $X=X_+\cup X_+^{-1}$. При этом полугруппа $X_+$
 индуцирует в $X$  линейный порядок,
согласованный со структурой
группы, по правилу $\xi\leq\chi:=\chi\xi^{-1}\in X_+$. Далее мы
положим $X_-:= X\setminus X_+=X_+^{-1}\setminus\{\textbf{1}\}$.

Как известно,  абелева группа $X$ может быть линейно
упорядочена тогда и только тогда, когда она не имеет кручения
(см., например, \cite{KK}), что, в свою очередь, равносильно
тому, что её группа характеров $G$ компактна и связна \cite{Pont} (при этом линейный порядок
в $X$, вообще говоря, не единственен). В приложениях в роли  $X$ часто выступают подгруппы аддитивной группы $\mathbb{R}^n$,
наделенные дискретной топологией, так что $G$ является боровской компактификацией группы  $X$.
В частности, в качестве  $X$ можно взять группу $\mathbb{Z}^n$, наделенную лексикографическим порядком. В этом случае $X$ имеет наименьший положительный элемент $(0,\dots,0,1)$, а $G=\mathbb{T}^n$ ---  $n$-мерный тор. Относительно других примеров см. \cite {SbMath}.

Через $\widehat{\varphi}$   мы будем
обозначать преобразование Фурье функции $\varphi$ из  $L^1(G)$, т.~е.
$$
\widehat{\varphi}(\xi)=\int\limits_G\varphi(x)\overline{\xi(x)}dx,\
\xi\in X.
$$

{\bf Определение 1.} {\it Пространство Харди} $H^p(G)\ (1\leq
p\leq\infty)$ {\it над} $G$ определяется следующим образом (см.,
например, \cite{Rud}):
$$
H^p(G)=\{f\in L^p(G):\widehat{f}(\chi)=0\ \forall\chi\in X_-\}.
$$
Обозначим через $H^2_-(G)$ ортогональное дополнение подпространства
$H^2(G)$ пространства $L^2(G)$. Тогда
$$H^2_-(G)=\{f\in L^2(G):\widehat{f}(\chi)=0\ \forall\chi\in X_+\}.$$
При этом $X$ является ортонормированным базисом пространства
$L^2(G)$,  $X_+$ -- ортонормированным базисом пространства
$H^2(G)$,  а  $X_-$ -- ортонормированным базисом пространства
$H^2_-(G)$. Через $P_+$ и $P_-$ мы будем обозначать ортопроекторы из
$L^2(G)$ на $H^2(G)$ и $H^2_-(G)$ соответственно.

Напомним определение преобразования Гильберта на группе
$G$. Мы ограничимся  случаем линейного порядка на $X$, принадлежащим
С. Бохнеру и Г. Хелсону (см., например, \cite[глава
8]{Rud}), более общая теория построена в \cite[глава
6]{garman} и \cite{ijpam}.  Для любой функции $u$ из $L^2(G,\mathbb{R})$
существует единственная функция $\widetilde{u}$ из
$L^2(G,\mathbb{R})$, такая, что $\widehat{\widetilde{u}}(1)=0$ и
$u+{\it i}\widetilde{u}\in H^2(G)$. Функция $\widetilde{u}$
называется {\it гармонически сопряженной с $u$}.
 Линейное отображение $\mathcal{H}$,
получаемое в результате продолжения отображения
$u\mapsto \widetilde{u}$ на (комплексное) $L^2(G)$ по
линейности, называется {\it преобразованием Гильберта} на группе
$G$. Этот оператор ограничен в $L^2(G)$.

Доказательство следующей леммы можно найти в \cite{garman} или \cite{ijpam}.

{\bf Лемма 1}. {\it Если $u\in L^2(G)$ ,
то $\widehat{\widetilde{u}}=- i{\rm sgn}_{X_+}\cdot\widehat{u}$, где
${\rm sgn}_{X_+}:=1_{X_+}-1_{X_+^{-1}}.$}

\begin{center}
{\bf  3. Некоторые свойства пространств функций
ограниченной средней осцилляции на группе $G$}
\end{center}

Следующие определения мотивированы
известной теоремой Ч. Феффермана \cite{F} (см. также \cite[с.
189]{Nik1}).

\textbf{Определение 2.} Определим пространства $BMO(G)$
\textit{функций ограниченной средней осцилляции} и $BMOA(G)$
\textit{функций ограниченной средней осцилляции аналитического типа}
на группе  $G$ следующим образом:
$$
BMO(G):=\{f+\widetilde{g}: f,g\in L^{\infty}(G)\},
$$
$$
BMOA(G):=BMO(G)\cap H^1(G),
$$
$$
\|\varphi\|_{BMO}:=\inf\{\|f\|_{\infty}+\|g\|_{\infty}:\varphi=f+\widetilde{g}, f,g\in L^{\infty}(G)\}
\ (\varphi\in BMO(G)).
$$
Норма в пространстве $BMOA(G)$
определяется так же, как и в $BMO(G)$.

  В идующем ниже предложении  вложения функциональных пространств понимаются в смысле \cite[с.
124]{BIN}.

{\bf Предложение 1. }{\it Имеют место следующие непрерывные вложения
пространств:
$$
L^{\infty}(G)\subset BMO(G)\subset
\bigcap\limits_{1<p<\infty}L^p(G).
$$}

{\it Доказательство.} Рассмотрим произвольную функцию $\varphi\in
BMO(G)$. По определению она представима в виде $\varphi=f+\widetilde{g}$, где
$f,g\in L^{\infty}(G)\subset L^p(G)$ при
всех $p\in(1;\infty)$.  Кроме того, для этих $p$
имеем $\mathcal{H}g=\widetilde{g}\in L^p(G)$,  и
оператор $\mathcal{H}:L^p(G)\rightarrow L^p(G)$ ограничен (теорема М.~Рисса; см., например, \cite{Rud}, теорема 8.7.2).  Тот факт,
что $L^{\infty}(G)$ содержится в $BMO(G)$, сразу вытекает из
определения 2.
Таким образом, $L^{\infty}(G)\subset BMO(G)\subset
\bigcap_{1<p<\infty}L^p(G)$.

Вложение
$L^{\infty}(G)\rightarrow BMO(G)$
непрерывно в силу очевидного
неравенства $\|f\|_{BMO}\leq\|f\|_{\infty}\ (f\in  L^{\infty}(G))$.
Рассмотрим вложение
$BMO(G)\rightarrow L^p(G).$
Для любого $p\in (1,\infty)$ в силу вышеупомянутой теоремы М. Рисса   для некоторой константы $K_p$ при всех $\varphi\in
BMO(G), \varphi=f+\widetilde{g}, f,g\in L^{\infty}(G)$ справедливо неравенство

$$
\|\varphi\|_p\leq \|f\|_p+\|\widetilde{g}\|_p
\leq\|f\|_{\infty}+K_p\|g\|_{\infty}\leq(1+K_p)(\|f\|_\infty)+\|g\|_{\infty}).
$$
 Переходя к инфимуму по всевозможным
разложениям функции $\varphi$, получим непрерывность вложения
$BMO(G)$ в $L^p(G)$.

{\bf Предложение 2. }{\it Пространства $BMO(G)$ и $BMOA(G)$ являются
банаховыми.}

{\it Доказательство.} Рассмотрим абсолютно сходящийся в $BMO(G)$ ряд

$$
\sum\limits_{n=0}^{\infty}h_n,\ h_n\in BMO(G).\eqno (3.1)
$$
Согласно определению нормы пространства $BMO(G)$, для любого $n$
существуют функции $f_n,g_n\in L^{\infty}$ такие, что
$h_n=f_n+\widetilde{g_n}$
и
$$
\|f_n\|_{\infty}+\|g_n\|_{\infty}<\|h_n\|_{BMO}+\frac{1}{2^n}.
$$
Поэтому числовой ряд
$\sum_{n=0}^{\infty}\left(\|f_n\|_{\infty}+\|g_n\|_{\infty}\right)$, а вслед за ним и ряды
$\sum_{n=0}^{\infty}\|f_n\|_{\infty}$ и
$\sum_{n=0}^{\infty}\|g_n\|_{\infty}$,  сходятся. Значит, в силу полноты пространства $L^{\infty}(G)$  ряды
$$
\sum\limits_{n=0}^{\infty}f_n\makebox{ и }\sum\limits_{n=0}^{\infty}g_n\eqno (3.2)
$$
сходятся в нем, т.~е. существуют функции $f,g\in L^{\infty}(G)$
такие, что
$$
\lim\limits_{N\rightarrow\infty}\left\|\sum\limits_{k=0}^Nf_k-f\right\|_{\infty}=0,\ \lim\limits_{N\rightarrow\infty}\left\|\sum\limits_{k=0}^Ng_k-g\right\|_{\infty}=0.
$$
Докажем, что  ряды (3.2) сходятся в пространстве $BMO(G)$.
По предложению 1  $L^{\infty}(G)\subset BMO(G)$, и
$\|f\|_{BMO}\leq\|f\|_{\infty}$ для любой функции $f\in
L^{\infty}(G)$. Поэтому
$$
\left\|\sum\limits_{k=0}^Nf_k-f\right\|_{BMO}\leq\left\|\sum\limits_{k=0}^Nf_k-f\right\|_{\infty}\rightarrow 0.
$$
Следовательно, $\sum_{k=0}^Nf_k\rightarrow f$ в пространстве
$BMO(G)$.
Аналогично для второго ряда. Далее, поскольку
$\mathcal{H}L^{\infty}(G)\subset BMO(G)$ и
$\|\mathcal{H}{g}\|_{BMO}\leq\|g\|_{\infty}$ для любой функции $g\in
L^{\infty}(G)$, то

$$
\left\|\sum\limits_{k=0}^N\widetilde{g_k}-\widetilde{g}\right\|_{BMO}=
\left\|\mathcal{H}\left(\sum\limits_{k=0}^N g_k-g\right)\right\|_{BMO}\leq\left\|\sum\limits_{k=0}^Ng_k-g\right\|_{\infty}\rightarrow 0,
$$
а значит $\sum_{k=0}^N\widetilde{g_k}\rightarrow
\widetilde{g}$ в пространстве $BMO(G)$.
Так как $\sum_{n=0}^{N}h_n=
\sum_{k=0}^Nf_k+\sum_{k=0}^N\widetilde{g_k}$, то отсюда следует сходимость ряда
(3.1) в этом пространстве. Согласно критерию полноты
нормированных пространств, пространство $BMO(G)$ является банаховым.

Теперь рассмотрим абсолютно сходящийся в $BMOA(G)$ ряд
$$
\sum\limits_{n=0}^{\infty}h_n,\ h_n\in BMOA(G).\eqno (3.3)
$$
Как показано выше, ряд (3.3) сходится к функции $h\in BMO(G)$.
Заметим, что в силу предложения 1 $BMO(G)\subset L^{2}(G)\subset L^1(G)$, причем ряд (3.3) сходится в $L^2(G)$.
Из определения
$BMOA(G)$ следует, что $\widehat{h_n}(\chi)=0$
для всех $\chi\in X_-,\ n\geq 0$. Тогда для всех $\chi\in X_-$ имеем по теореме Планшереля
$
\widehat{h}(\chi)=\sum_{n=0}^{\infty}\widehat{h_n}(\chi)
=0,$ а потому  $h\in BMOA(G)$. Следовательно, $BMOA(G)$
банахово, что и требовалось доказать.

{\bf Лемма 2.} \textit{Для любой функции $\psi$ из $L^{2}(G)$ верно
равенство}
$$
\widetilde{\psi}=-i(P_+\psi-P_-\psi-\widehat{\psi}(\textbf{1})).
$$

Доказательство. По теореме Планшереля $\psi=\sum_{\chi\in
X}\widehat{\psi}(\chi)\chi$ (ряд сходится в $L^2(G)$). Положим
$\varphi:=- i (P_+\psi-P_-\psi-\widehat{\psi}(\textbf{1}))$. Тогда имеем для любого
$\xi\in X$  с учетом ортогональности характеров
$$
\widehat{\varphi}(\xi)=\int\limits_G\varphi(x)\overline{\xi(x)}dx=-
i\sum\limits_{\chi\in
X_+}\widehat{\psi}(\chi)\int\limits_G\chi(x)\overline{\xi(x)}dx
$$
$$
+ i\sum\limits_{\chi\in
X_-}\widehat{\psi}(\chi)\int\limits_G\chi(x)\overline{\xi(x)}dx+
i\widehat{\psi}(\textbf{1})\int\limits_G \overline{\xi(x)}dx
$$
$$
=- i1_{X_+}(\xi)\widehat{\psi}(\xi)+ i1_{X_-\setminus{\{
\textbf{1}\}}}(\xi)\widehat{\psi}(\xi)+ i1_{\{ \textbf{1}\}}(\xi)\widehat{\psi}(\xi)=-
i{\rm sgn}_{X_+}\widehat{\psi}(\xi).
$$
Так как $\widehat{\widetilde{\psi}}=- i{\rm
sgn}_{X_+}\widehat{\psi}$  по лемме 1,  то
$\widehat{\widetilde{\psi}}(\xi)=\widehat{\varphi}(\xi)$ для всех
$\xi\in X$, и значит $\widetilde{\psi}=\varphi$. Лемма доказана.

{\bf Лемма 3. }{\it Имеет место равенство}
$$
BMOA(G)=BMO(G)\cap H^2(G).
$$

{\it Доказательство.} Так как $G$ компактно, то $L^2(G)\subset
L^1(G)$, а значит $H^2(G)\subset H^1(G)$. Поэтому
$BMO(G)\cap H^2(G)\subset BMOA(G).$
С другой стороны, если $\varphi\in BMOA(G)$, то
$\varphi\in BMO(G)\subset L^2(G)$ и $\widehat{\varphi}(\chi)=0$ для всех
$\chi\in X_-.$ Значит, $\varphi\in BMO(G)\cap H^2(G)$. Лемма доказана.

{\bf Предложение 3. }{\it Справедливы следующие утверждения:

$
(1)\  BMO(G)=P_-L^{\infty}(G)+P_+L^{\infty}(G),
$
причем
$$
\|\varphi\|_{\ast BMO}:=\inf\{\max(\|f_1\|_\infty, \|g_1\|_\infty):\varphi=P_-f_1+P_+g_1, f_1, g_1\in L^\infty(G)\} -
$$
эквивалентная норма;

$
(2)\  BMOA(G)=P_+L^{\infty}(G),
$
причем
$$
\|\varphi\|_{\ast BMOA}:=\inf\{\|g_1\|_\infty:\varphi=P_+g_1,\ g_1\in L^\infty(G)\} -
$$
эквивалентная норма.}

{\it Доказательство.} (1) Рассмотрим произвольную функцию $\varphi$ из
$BMO(G)$, $\varphi=f+\widetilde{g}, f,g\in
L^{\infty}(G)$. В силу леммы 2
$$
\varphi=P_+f+P_-f-{\it i}(P_+g-P_-g-\widehat{g}({\bf 1}))=
$$
$$
=P_+(f-{\it i}g+{\it i}\widehat{g}({\bf 1}))+P_-(f+{\it i}g)=P_+g_1+P_-f_1
$$
где функции $g_1=f-{\it i}g+{\it i}\widehat{g}({\bf 1})\makebox{ и
}f_1=f+{\it i}g \makebox{ принадлежат } L^{\infty}(G).$
Поскольку $\|f_1\|_\infty\leq \|f\|_\infty+\|g\|_\infty, \|g_1\|_\infty\leq 2(\|f\|_\infty+\|g\|_\infty)$   ,  то  $\|\varphi\|_
{\ast BMO}\leq 2\|\varphi\|_
{BMO}$.

Обратно, рассмотрим произвольную функцию
$$
\varphi\in P_+L^{\infty}(G)+P_-L^{\infty}(G),\varphi=P_-f_1+P_+g_1.
$$
Из леммы 2 следует, что
$$
P_-h=\frac{1}{2}(h-{\it i}\widetilde{h}-\widehat{h}({\bf 1}))\makebox{ и }P_+h=\frac{1}{2}(h+{\it i}\widetilde{h}+\widehat{h}({\bf 1}))
$$
для всех $h\in L^{\infty}(G)$. Поэтому
$$
\varphi=P_-f_1+P_+g_1=\frac{1}{2}(f_1-{\it i}\widetilde{f_1}-\widehat{f_1}({\bf 1}))+\frac{1}{2}(g_1+{\it i}\widetilde{g_1}+\widehat{g_1}({\bf 1}))=
$$
$$
=\frac{1}{2}(f_1+g_1-\widehat{f_1}({\bf 1})+\widehat{g_1}({\bf 1}))+\frac{1}{2}({\it i}\widetilde{g_1}-{\it i}\widetilde{f_1})=f+\widetilde{g}\in BMO(G),
$$
где $f=\frac{1}{2}(f_1+g_1-\widehat{f_1}({\bf 1})+\widehat{g_1}({\bf
1}))\makebox{ и }g=\frac{1}{2}({\it i}g_1-{\it i}f_1),f,g\in
L^{\infty}(G).$

Далее,
$$
\|\varphi\|_
{BMO}\leq \|f\|_\infty+\|g\|_\infty\leq \frac{3}{2}(\|f_1\|_\infty+\|g_1\|_\infty)\leq 3\|\varphi\|_
{\ast BMO}.
$$
Таким образом, окончательно имеем
$$
\frac{1}{3}\|\varphi\|_
{BMO}\leq \|\varphi\|_
{\ast BMO}\leq 2\|\varphi\|_
{BMO},
$$
и первое утверждение доказано.

(2) Второе равенство следует из первого и леммы 3:
$$
BMOA(G)=(P_-L^{\infty}(G)+P_+L^{\infty}(G))\cap H^2(G)=P_+L^{\infty}(G).
$$
 А так как $\|\varphi\|_
{\ast BMO}=\|\varphi\|_
{\ast BMOA}$ при $\varphi\in BMOA(G)$, утверждение о нормах следует из (1) (впрочем, как и в случае (1), легко проверить напрямую, что
$$
\frac{2}{3}\|\varphi\|_
{BMO}\leq \|\varphi\|_
{\ast BMOA}\leq 2\|\varphi\|_
{BMO}).
$$
Предложение  доказано.

{\bf Следствие 1.}\textit{ Справедливо равенство $BMOA(G)=P_+BMO(G)$.}

{\bf Лемма 4. }{\it Если $\varphi\in BMO(G)$, то и
$\overline{\varphi}\in BMO(G)$.}

{\it Доказательство.} Пусть $\varphi\in BMO(G)$. Тогда по предложению 3
для некоторых $f,g\in L^{\infty}(G)$ имеем
$$\varphi=P_+g+P_-f=\sum\limits_{\chi\in X_+}a_{\chi}\chi+\sum\limits_{\xi\in X_-}b_{\xi}\xi,$$
где $f=\sum_{\xi\in X}b_{\xi}\xi, g=\sum_{\chi\in
X}a_{\chi}\chi$ -- разложения Фурье.
Поэтому, снова воспользовавшись предложением 3, имеем:
$$
\overline{\varphi}=\sum\limits_{\chi\in X_+}\overline{a_{\chi}}\bar{\chi}+\sum\limits_{\xi\in X_-}\overline{b_{\xi}}\bar{\xi}
=\sum\limits_{\chi\in X_+\setminus\{{\bf
1}\}}\overline{a_{\chi}}\bar{\chi}+\overline{a_{1}}{\bf
1}+\sum\limits_{\xi\in X_-}\overline{b_{\xi}}\bar{\xi}=
$$
$$
=P_-f_1+P_+g_1\in BMO(G),
$$
где $f_1=\sum_{\xi\in
X}\overline{a_{\xi}}\xi,g_1=\overline{a_{1}}{\bf
1}+\sum_{\chi\in X\backslash\{{\bf
1}\}}\overline{b_{\chi}}\chi, f_1,g_1\in L^{\infty}(G)$. Лемма
доказана.

\begin{center}
{\bf  4. Связь с операторами Ганкеля}
\end{center}

{\bf Определение 3.} Пусть $k$ --- функция на $X_+$. \textit{Ганкелевой формой на $X_+$ c ядром $k$} называют комплексную
билинейную форму вида
$$
A(a,b)=\sum\limits_{\chi,\eta\in X_+}k(\chi\eta)a(\chi)b(\eta),
$$
определенную первоначально на финитных  функциях  на $X_+$.

Далее важную роль будет играть следующий результат, обобщающий классическую теорему Нехари \cite{Neh}.

{\bf Теорема Нехари-Вонга} \cite{Wong}. {\it Ганкелева форма с ядром $k$ на
$X_+$ ограничена тогда и только тогда, когда ее ядро имеет вид
$k(\chi)=\widehat{\varphi}(\overline\chi),\ \chi\in X_+,$
 где $\varphi\in L^{\infty}(G)$.  При этом
$\|\varphi\|_{\infty}\leq M,$  где $M$ --- константа ограниченности
формы $A.$ Кроме того, норма формы $A$ равняется
$\|\varphi\|_{\infty}$ для некоторой функции $\varphi$,
удовлетворяющей указанным выше условиям}.

Нам понадобятся две реализации операторов Ганкеля.

{\bf Определение 4.}  Оператор $\Gamma:l_2(X_+)\to l_2(X_+)$,
определенный первоначально на финитных  функциях  на $X_+$,
называется {\it ганкелевым (оператором Ганкеля) в} $l_2(X_+)$, если
существует функция $a=a_\Gamma$ на $X_+$ такая, что для всех
$\chi,\xi\in X_+$ выполняется равенство
$$
\langle\Gamma 1_{\{\chi\}},1_{\{\xi\}}\rangle=a(\chi\xi)
$$
($1_{\{\xi\}}$ --- индикатор одноточечного подмножества $\{\xi\}\subset X$; угловые скобки обозначают скалярное произведение в $l_2(X_+)$).

{\bf Определение 5.} Пусть $\varphi\in L^{2}(G)$. {\it Оператором
Ганкеля (ганкелевым оператором) в $H^2(G)$ с символом $\varphi$} назовем
оператор $H_{\varphi}:H^2(G)\rightarrow H^2_-(G)$, определяемый
первоначально на  пространстве  тригонометрических
полиномов аналитического типа (линейных комбинациях характеров из
$X_+$) равенством
$$
H_{\varphi}f=P_-({\varphi}f).
$$

{\bf Лемма 5.}\textit{ Пусть $X$ содержит наименьший положительный
элемент $\chi_1$.}

(а) \textit{Справедливо равенство $X_-=X_+^{-1}\chi_1^{-1}$.
Следовательно, отображения
$$
i:\{1_{\{\chi\}}\}_{\chi\in X_+}\rightarrow
X_+:1_{\{\chi\}}\mapsto\chi,\ j:\{1_{\{\chi\}}\}_{\chi\in
X_+}\rightarrow X_-:1_{\{\chi\}}\mapsto{\chi}^{-1}\chi_1^{-1}$$
 продолжаются единственным образом  до изоморфизмов гильбертовых пространств}
$$
i:l_2(X_+)\rightarrow H^2(G),\ j:l_2(X_+)\rightarrow H_-^2(G).
$$

(б) \textit{Оператор $H_\varphi$ с символом $\varphi\in L^\infty (G)$
унитарно эквивалентен ограниченному  ганкелеву оператору в
$l_2(X_+)$ и обратно.}

Доказательство. (а) Так как $X_+\setminus\chi_1 X_+=\{\textbf{1}\}$, то
$\chi_1 X_+=X_+\setminus\{\textbf{1}\}$;  первое утверждение получается
отсюда переходом к обратным элементам. Поэтому,
 когда $\chi$ пробегает  $X_+$,\ $\chi^{-1}\chi_1^{-1}$  пробегает  $X_-$, а потому $j$ (и, очевидно,  $i$) являются биекциями. Второе утверждение  теперь следует из  теоремы Рисса-Фишера.

(б)  Положим $\Gamma=j^{-1}H_{\varphi}i$. Тогда этот оператор
является ганкелевым в $l_2(X_+)$. Действительно, для любых
$\chi,\xi\in X_+$
$$
\langle\Gamma1_{\{\chi\}},1_{\{\xi\}}\rangle_{l_2}= \langle
j^{-1}P_-(\varphi\chi),j^{-1}(\overline{\xi}\chi_1^{-1})\rangle_{l_2}=\langle
\varphi\chi,P_-(\overline{\xi}\chi_1^{-1})\rangle_{L^2}
$$
$$
=\langle \varphi\chi,\overline{\xi}\chi_1^{-1}\rangle_{L^2}=
\widehat{\varphi}(\overline{\chi\xi}\chi_1^{-1})=a_\Gamma(\chi\xi),
$$
где $a_\Gamma(\chi):=\widehat{\varphi}(\overline{\chi}\chi_1^{-1}),\
\chi\in X_+$.

Если еще $\varphi\in L^\infty (G)$, то оператор $\Gamma$ ограничен
вслед за $H_\varphi$.

Обратно, если   $\Gamma$ есть ограниченный ганкелев оператор в
$l_2(X_+)$, то положим $A=j \Gamma i^{-1}$. Тогда $ \langle
A\chi,\bar\xi\rangle=\langle\Gamma 1_{\{\chi\}},
j^{-1}\xi^{-1}\rangle=\langle\Gamma 1_{\{\chi\}},
1_{\{\zeta\}}\rangle,$ где $\chi\in X_+, \xi\in X_+\backslash\{
\textbf{1}\},\ \xi=\zeta\chi_1,\ \zeta\in X_+$. Следовательно, $\langle
A\chi,\bar\xi\rangle=a(\chi\zeta)=a(\chi\xi\chi_1^{-1})$ зависит
лишь от $\chi\xi$, а потому, снова  в силу следствия теоремы 2.1 из
\cite{Dyba},  оператор $A$ имеет вид $H_{\varphi}$, где  $\varphi\in
L^\infty(G)$. Лемма доказана.

{\bf Лемма 6.} {\it Оператор Ганкеля $\Gamma$ в
$l_2(X_+)$  ограничен тогда и только тогда, когда существует функция $\psi\in
L^{\infty}(G)$ такая, что $a_{\Gamma}(\chi)=\widehat{\psi}(\chi)$
для любого $\chi\in X_+$.}

{\it Доказательство.} Рассмотрим билинейную форму ($f,g$ --- $\mbox{
финитные  функции  на } X_+$)
$$
A(f,g):=\langle\Gamma f;\bar{g}\rangle.
$$
Поскольку $f=\sum_{\chi\in X_+}f(\chi)1_{\{\chi\}},g=\sum_{\xi\in
X_+}g(\xi)1_{\{\xi\}}$, то
$$
\langle\Gamma f,\bar{g}\rangle= \sum\limits_{\chi,\xi\in
X_+}a_\Gamma(\chi\xi)f(\chi)g(\xi),
$$
а значит форма $A$ является ганкелевой. По теореме  Нехари-Вонга
форма $A$, а вместе с ней и оператор $\Gamma$, ограничена тогда и
только тогда, когда существует функция $\psi_1\in L^{\infty}(G)$
такая,что для любого $\chi\in X_+$ выполняется равенство
$a_{\Gamma}(\chi)=\widehat{\psi_1}(\bar{\chi})$, и осталось положить в этой формуле $\psi(x)=\psi_1(x^{-1})$.  Лемма доказана.

{\bf Теорема 1.} {\it Пусть в группе $X$ существует наименьший
положительный элемент. Имеют место следующие равенства:
$$
BMOA(G)=\{\varphi\in H^2(G):H_{\overline{\varphi}}\makebox{ ограничен}\},\eqno(4.1)
$$
$$
BMO(G)=\{\varphi\in L^2(G):H_{\overline{\varphi}},H_{\varphi}\makebox{ ограничены}\}. \eqno(4.2)
$$

Кроме того, формулы
$$
\|\varphi\|_H:=\|H_{\overline \varphi}\|+\|H_{\varphi}\|,\ \|\varphi\|^\prime_H:=\|H_{\overline \varphi}\|
$$
задают полунормы на пространствах $BMO(G)$ и $BMOA(G)$ соответственно,  нулевые подпространства которых совпадают с  пространством $\mathbb{C}\cdot\textbf{1}$ постоянных функций,  причем нормы, отвечающие этим полунормам после соответствующей факторизации, эквивалентны  норме $\|\cdot\|_{BMO}$.}

{\it Доказательство} разобьем на несколько шагов.

1. Докажем сначала, что \textit{для функции $\varphi\in L^2(G)$ оператор $H_{\varphi}$ ограничен тогда и только
тогда, когда $P_-\varphi\in BMO(G)$}. Для этого воспользуемся
обозначениями и результатами, содержащимися в  лемме 5 и ее
доказательстве. Сперва покажем, что $H_{\varphi}$ ограничен тогда и
только тогда , когда
$$\mbox{\it существует функция } \psi_1\in
L^{\infty}(G) \mbox{\it\ такая, что }
\widehat{\psi_1}|X_-=\widehat{\varphi}|X_-.\eqno (*)$$
В самом деле,  оператор
$H_{\varphi}$ ограничен тогда и только тогда, когда ограничен
оператор $\Gamma=j^{-1}H_{\varphi}i$. В доказательстве утверждения
(б) леммы 5 показано, что последний оператор является ганкелевым в
$l_2(X_+)$. Значит, по лемме 6 оператор $\Gamma$ ограничен тогда и
только тогда, когда существует функция $\psi\in L^{\infty}(G)$
такая, что
$\widehat{\psi}(\chi)=a_\Gamma(\chi)=\widehat{\varphi}(\overline{\chi}\chi_1^{-1})$
для всех $\chi\in X_+$  (последнее равенство установлено в
доказательстве леммы 5). Положим
$\psi_1(x)=\chi_1^{-1}(x)\psi(x^{-1})$. Тогда
$\widehat{\psi_1}(\bar{\chi}\chi_1^{-1})=\widehat{\psi}(\chi)$ для
любого $\chi\in X_+$, а потому
$\widehat{\psi_1}(\bar{\chi}\chi_1^{-1})=\widehat{\varphi}(\bar{\chi}\chi_1^{-1})$.
Таким образом, равенство
$\widehat{\psi}(\chi)=\widehat{\varphi}(\overline{\chi}\chi_1^{-1})$
для всех $\chi\in X_+$ эквивалентно равенству
$\widehat\psi_1(\xi)=\widehat{\varphi}(\xi)$ для всех $\xi\in X_-$,
и равносильность ограниченности оператора $H_{\varphi}$ и условия
$(*)$ доказана.

Докажем теперь, что $P_-\varphi\in BMO(G)$ тогда и только тогда,
когда выполнено условие $(*)$.

Пусть выполнено условие $(*)$. Из леммы 2 следует, что
$$P_-\psi_1=(\psi_1-\widehat{\psi_1}(\textbf{1}))/2+
(\widetilde{i\psi_1/2})\in BMO(G).$$
Это влечет необходимость,
поскольку
$$
P_-\psi_1=\sum\limits_{\xi\in
X_-}\widehat{\psi_1}(\xi)\xi=\sum\limits_{\xi\in
X_-}\widehat{\varphi}(\xi)\xi=P_-\varphi.
$$

Пусть теперь $P_-\varphi\in BMO(G)$. Тогда
$P_-\varphi=f+\widetilde{g}$, где $f, g\in L^\infty(G)$. Снова
применяя лемму 2, имеем
$$
P_-\varphi=P_-f+P_-(- i(P_+g-P_-g-\widehat{g}(\textbf{1})))= P_-(f+ ig+
i\widehat{g}(\textbf{1})).
$$
Рассмотрим функцию $\psi_1=f+ ig+ i\widehat{g}(\textbf{1})\in L^{\infty}(G)$.
Так как $P_-\psi_1=P_-\varphi$,  то
$\widehat{\psi_1}(\xi)=\widehat{\varphi}(\xi)$ для всех $\xi\in
X_-$, т.~е. условие $(*)$ выполняется.

2. Докажем  равенство (4.1). Пусть $\varphi\in BMOA(G)$. Из лемм 3  и 4 следует, что
$\varphi\in H^2(G)$, причем $\overline{\varphi}\in
BMO(G)$. Тогда по предложению  3 $\overline{\varphi}=P_+f+P_-g$ для
некоторых $f,g\in L^{\infty}(G)$, а значит
$P_-\overline{\varphi}=P_-g\in BMO(G)$.
Таким образом, по доказанному на шаге 1 оператор $H_{\overline{\varphi}}$ ограничен, и, стало быть,
$$
BMOA(G)\subset \{\varphi\in H^2(G):H_{\overline{\varphi}}\makebox{
ограничен}\}.
$$

Обратно, рассмотрим произвольную функцию $\varphi\in H^2(G)$, для которой оператор $H_{\overline{\varphi}}$ ограничен. По доказанному на шаге 1
$P_-{\overline{\varphi}}\in BMO(G)$. Поскольку $\varphi\in H^2(G)$, то
$
\varphi=\sum_{\chi\in X_+}a_{\chi}\chi.
$
Тогда
$$
\overline{\varphi}=\sum\limits_{\chi\in X_+}\overline{a_{\chi}}\bar{\chi}=\overline{a_1}{\bf 1}+P_-\overline{\varphi}\in BMO(G).
$$
По лемме 4 $\varphi\in BMO(G)$, что с учетом леммы 3 завершает
доказательство  равенства (4.1).

3.  Пусть теперь $\varphi\in BMO(G)$. По предложению 1
$\varphi\in L^2(G)$. В то же время из предложения 3 следует, что
$$
BMO(G)=BMOA(G)+\overline{BMOA(G)}.
$$
Если $\varphi=f_1+\overline{f_2}, f_1,f_2\in BMOA(G)$, то $H_{\overline{\varphi}}=H_{\overline{f_1}}+H_{f_2}=H_{\overline{f_1}}$ и
аналогично $H_{\varphi}=H_{\overline{f_2}}$, причем из
равенства (4.1) следует, что эти операторы ограничены.
Следовательно,
$$
BMO(G)\subset\{\varphi\in L^2(G):H_{\overline{\varphi}},H_{\varphi}\makebox{ ограничены}\}
$$

Обратно, пусть $\varphi\in L^2(G)$ и операторы
$H_{\overline{\varphi}},H_{\varphi}$ ограничены. Так как
 $\varphi=f_1+f_2'$, где $f_1\in
H^2(G)$ и $f_2'\in H^2_-(G)$, то
$H_{\varphi}=H_{f_1}+H_{f_2'}=H_{\overline{f_2}}$, где
$f_2=\overline{f_2'}\in H^2(G)$.
 Так как оператор
$H_{\varphi}$ ограничен, то, снова применяя   (4.1), получаем, что
$f_2\in BMOA(G)$. Аналогично доказывается, что $f_1\in
BMOA(G)$. Следовательно, $\varphi\in BMOA(G)+\overline{BMOA(G)}=BMO(G)$, и  равенство (4.2) доказано.

4. Поскольку $\|\varphi\|_H=\|\varphi\|^\prime_H$ при $\varphi\in BMOA(G)$, достаточно доказать эквивалентность норм $\|\cdot\|_H$ и $\|\cdot\|_{BMO}$. Ясно, что $\|\cdot\|_H$ является полунормой. Далее,
если $\|\varphi\|_H=0$, то $P_-\overline{\varphi}=P_-{\varphi}=0$, что влечет $\varphi=const$. Обратное утверждение очевидно.
Ниже для $\varphi$ из $BMO(G)$  через
$\dot \varphi:=\varphi+\mathbb{C}\cdot\textbf{1}$ мы обозначаем элемент факторпространства $BMO(G)/\mathbb{C}\cdot\textbf{1}$.
 Пусть $\varphi=P_-f+P_+g$, где $f, g\in L^\infty(G)$ (предложение 3). Используя разложения в ряд Фурье по характерам в $L^2(G)$, легко проверить, что
 $$
  \overline{\varphi}=P_+\overline{f}+P_-\overline{g}+(\widehat{\overline{g}}(\textbf{1})-\widehat{\overline{f}}(\textbf{1}))\textbf{1}.
 $$
Следовательно, $P_-{\overline \varphi}=P_-{\overline g}$. Поэтому
$$
\|\varphi\|_H=\|H_{P_-\overline \varphi}\|+\|H_{P_- \varphi}\|=\|H_{P_-\overline g}\|+\|H_{P_-f}\|=
$$
$$=\|H_{\overline g}\|+\|H_{f}\|\leq \|g\|_\infty+\|f\|_\infty\leq 2\max(\|g\|_\infty,\|f\|_\infty).
$$
Переходя к инфимуму по $f, g$, получаем с учетом предложения 3, что  $\|\varphi\|_H\leq 2\|\varphi\|_{\ast BMO}\leq 4\|\varphi\|_{BMO}$. Отсюда следует, что для любого комплексного $c$ имеем
$
\|\dot\varphi\|_H=\|\varphi+c\cdot\textbf{1}\|_H\leq 4\|\varphi+c\cdot\textbf{1}\|_{BMO},
$
а, стало быть,
$$
\|\dot \varphi\|_H\leq 4\inf_{c\in \mathbb{C}} \|\varphi+c\cdot\textbf{1}\|_{BMO}=4\|\dot\varphi\|_{BMO}.
$$
Ввиду предложения 2 для доказательства эквивалентности норм теперь достаточно показать, что пространство $(BMO(G),\|\cdot\|_H)$ (а потому и $(BMO(G)/\mathbb{C}\cdot\textbf{1},\|\cdot\|_H)$) банахово. Для этого рассмотрим фундаментальную последовательность $(\varphi_n)\subset (BMO(G),\|\cdot\|_H)$. Последовательность операторов $(H_{\overline\varphi_n})$ равномерно сходится к некоторому ограниченному оператору $A:H^2(G)\to H^2_-(G)$, и при этом последовательность операторов $(H_{\varphi_n})$ равномерно сходится к некоторому ограниченному оператору $B:H^2(G)\to H^2_-(G)$. Как отмечено выше, если
$\varphi_n=P_-f_n+P_+g_n,\ f_n,g_n\in L^\infty(G)$, то  $P_-\overline{\varphi_n}=P_-\overline g_n$, а потому  $H_{\overline\varphi_n}=H_{\overline g_n}$. В силу обобщения теоремы Нехари из \cite{Dyba}, для любого $\chi\in X_+$ выполняется равенство
$H_{\overline g_n}S_\chi=P_-S_\chi H_{\overline g_n}$
($S_\chi:L^2(G)\to L^2(G), f\mapsto \chi f$). В пределе получаем $AS_\chi=P_-S_\chi A$ для любого $\chi\in X_+$, откуда, снова в силу  обобщенной теоремы Нехари, $A=H_g$ для некоторой функции $g\in L^\infty(G)$. Аналогично $B=H_{\overline g_1}$ для некоторой функции $g_1\in L^\infty(G)$. Непосредственно проверяется, что для сопряженных операторов справедливо равенство $H_{\varphi_n}^*=P_+H_{\overline\varphi_n}|H^2_-(G)$. Отсюда следует, что $B^*=P_+A|H^2_-(G)$, т.~е. $H_{\overline{g_1}}^*=P_+H_{g}|H^2_-(G)=H_{\overline{g}}^*$. Значит, $H_{\overline{g}}=H_{\overline{g_1}}=B$.
 При этом
$$
\|\varphi_n-\overline{g}\|_H=\|H_{\overline\varphi_n}-H_{g}\|+\|H_{\varphi_n}-H_{\overline g}\|=\|H_{\overline\varphi_n}-A\|+\|H_{\varphi_n}-B\|\to 0
$$
$(n\to\infty)$, что завершает доказательство эквивалентности норм для $BMO(G)$.
Теорема доказана.

\textbf{Следствие 2.} \textit{Пусть в группе $X$ существует наименьший
положительный элемент,   $\varphi\in L^2(G)$. Оператор $H_\varphi$ ограничен тогда и только тогда, когда $H_\varphi=H_g$ для некоторой функции $g\in L^\infty(G)$.}

\textit{Доказательство. } Как показано на шаге 1 доказательства теоремы 1, оператор $H_\varphi$ ограничен тогда и только тогда, когда $P_-\varphi\in BMO(G)$. В силу предложения 3 $P_-\varphi=P_-g$ для некоторой функции $g\in L^\infty(G)$, а потому $H_\varphi=H_g$.
Обратное утверждение очевидно.

\newpage
{\footnotesize Рассматривается  обобщение понятия функции ограниченной средней осцилляции и ганкелева оператора на случай компактных абелевых групп с линейно упорядоченной  группой характеров.  Дается описание пространств функций ограниченной средней осцилляции и  функций ограниченной средней осцилляции аналитического типа на таких группах в терминах ограниченности соответствующих операторов Ганкеля в предположении, что группа характеров содержит наименьший положительный элемент.}

{\footnotesize Ключевые слова:  оператор Ганкеля, ограниченный оператор,  ограниченная средняя осцилляция, линейно упорядоченная абелева группа, компактная абелева группа}.\\

 {\footnotesize R.V.~Dyba, A.R.~Mirotin. Functions of  bounded mean oscillation and Hankel operators  on compact abelian groups.

 Generalization of  functions of bounded mean oscillation and Hankel operators to the case of compact abelian groups with linearly ordered dual is considered. Spaces of functions of bounded mean oscillation and of  bounded mean oscillation of analytic type on such groups are described in terms  of boundedness of corresponding Hankel operators under the assumption that the dual group contains a minimal positive element.}

 {\footnotesize Key words: Hankel operator, bounded operator, bounded mean oscillation, linearly ordered abelian group, compact abelian group}.

\end{document}